\documentclass[12pt,a4paper,reqno]{amsart}    

\usepackage[a4paper,margin=1.5in, total={6in,9in}]{geometry}

\usepackage{graphicx}
\usepackage{pdflscape}
\usepackage{charter}

\usepackage{amsmath,bbm}
\usepackage{amssymb}
\usepackage{graphicx}
\usepackage{color}
\usepackage{latexsym}

\usepackage[english]{babel}
\usepackage[utf8]{inputenc}
\usepackage[T1]{fontenc}

\usepackage[plainpages=false,pdfpagelabels,colorlinks=true,citecolor=blue,hypertexnames=false]{hyperref}

\usepackage[x11names]{xcolor}
\usepackage{tcolorbox}
\usepackage{bm}

\catcode`\@=11
\def\section{\@startsection{section}{1}%
 \z@{.7\linespacing\@plus\linespacing}{.5\linespacing}%
 {\normalfont\bfseries\scshape\centering}}

\def\subsection{\@startsection{subsection}{2}%
 \z@{.5\linespacing\@plus\linespacing}{.5\linespacing}%
 {\normalfont\bfseries\scshape}}

\def\subsubsection{\@startsection{subsubsection}{3}%
 \z@{.5\linespacing\@plus\linespacing}{-.5em}
 {\normalfont\bfseries\itshape}}
\catcode`\@=12

\newtheorem{thm}{Theorem}
\newtheorem{fact}{Fact}

\newtheorem{example}{Example}
\newtheorem{prb}[thm]{Problem}
\newtheorem{definition}{Definition}

\def\Q{\mathbb{Q}}
\def\Sdup{\mathcal{S}_{\textsf{dup}}}

\def\opname#1{\operatorname{#1}}

\def\K{\mathbb{K}}

\def\code#1{\begin{tcolorbox}[colback=orange!5!white,colframe=orange!75!orange,]
    \begin{center} #1 \end{center} \end{tcolorbox}}
\def\mb#1{\mathbf{#1}}
\def\fracpuiseux{\overline{\Q}[[t^{\frac{1}{\star}}]]}
\def\algebraicclosure{\overline{\Q(t)}}
\def\ddesolver{\href{https://github.com/HNotarantonio/ddesolver}{DDE-Solver}\,}
\def\opname#1{\operatorname{#1}}
\def\LT{\opname{LT}}
\def\ann{$\operatorname{\textbf{annihilating\_polynomial}}$\,}

\title{\href{https://github.com/HNotarantonio/ddesolver}{DDE-Solver}: a \emph{Maple} package for Discrete Differential Equations}
\date{January 2025}

\begin{document}
\author[H. Notarantonio]{Hadrien Notarantonio}
  \address{H. Notarantonio: CNRS, IRIF, Université Paris Cité, 75013 Paris, France}
  \email{hadrien.notarantonio@irif.fr}
\date{}
\maketitle

\begin{abstract}
We introduce \href{https://mathexp.eu/notarantonio/}{DDE-Solver}, a
  Maple package designed for solving \emph{Discrete Differential Equations}~(DDEs).
  These equations are functional equations relating algebraically
  a formal power series~$F(t, u)$ with polynomial
  coefficients in a ``catalytic'' variable~$u$,
  with specializations of it with respect to the catalytic variable.
  Such equations appear in enumerative combinatorics, for instance
  in the enumeration of maps.
  Bousquet-Mélou and Jehanne showed in~$2006$ that
  when these equations are of a fixed point type in~$F$, then~$F$ is an algebraic
  series. In the same paper, they proposed
  a systematic method for computing annihilating polynomials of these series.
  Bostan, Safey El Din and the author of this paper
  recently designed new efficient algorithms
  for computing these witnesses~of algebraicity. This paper provides combinatorialists
  an automated tool in hand that solves DDEs using these algorithms. 
  We also compare the timings of all these algorithms on DDEs from the literature.
\end{abstract}

\section{Introduction}\label{sec:introduction}

Sequences of non-negative integers are ubiquitous in enumerative combinatorics.
For instance, when studying bicoloured maps (the faces are either coloured in purple or white)
such that the degree of each purple face
is $3$ and the degree of each white face is a multiple of $3$, the introduction
of the sequence $(c_n)_{n\in\mathbb{N}}
$ of such maps, called~$3$-constellations, with n purple faces starts with the numbers~\[1, 1, 6, 54, 594, 7371, 99144, \cdots.\]

For such a sequence, 
many questions arise:
does there exist a finite way of representing the infinite
amount of data given by this sequence?
Does there exist a closed formula for~$(c_n)_{n\in\mathbb{N}}$?
What is the value of~$c_{N}$, for~$N\in\mathbb{Z}_{\geq 0}$ large (e.g.~$N=3\cdot10^6$),
and how fast can it be computed? 
What is the asymptotic growth of~$c_n$ when $n\rightarrow\infty$? \\

\indent In order to answer such questions, a common method is to introduce and study
the properties of the associated generating
function. For instance for the above maps enumeration problem, one would
introduce~$G(t) := \sum_{n\in\mathbb{N}}c_nt^n \in\Q[[t]]$.
For sophisticated enumerations, it is usually hard to study directly the generating
series. Refining the initial enumeration of interest usually leads to introducing
a new variable called \emph{catalytic}, which leads to consider
a bivariate generating function, say~$F(t, u)$ for~$u$ the catalytic variable.
For well-chosen refinements, it is then possible to
write a functional equation relating algebraically~$F(t, u)$ and~$G(t)$.\\

\indent For the above toy example of bicoloured maps, one would typically refine the enumeration with
the sequence~$(c_{n, d})_{n, d \in\mathbb{N}}$, where~$c_{n, d}$
is number of 3-constellations having $n$ purple faces
and outer degree $3d$. Straightforwardly, one would also introduce its associated generating
function~$F(t, u) = \sum_{n, d\geq 0} c_{n, d}u^dt^n\in\Q[u][[t]]$. 
As~$\sum_{d=0}^\infty c_{n, d} = c_n$, it holds that~$F(t, 1) = G(t)$. An analysis of the construction
of such maps~\cite[Fig.~$7$]{BMJ06} yields the functional equation
\begin{align}\label{3const}
 F(t, u) = 1 + tuF(t, u)^3 + tu(2F(t, u)&+F(t, 1))\frac{F(t, u)-F(t, 1)}{u-1}\\\nonumber &+
  tu\frac{F(t, u) - F(t, 1) - (u-1)\cdot \partial_u F(t, 1)}{(u-1)^2}.
\end{align}
\indent
Note that Eq.~\eqref{3const} can be put in the form
\begin{equation}\label{DDE}
F(t, u) = f(u) + t \cdot Q(F(t, u), \Delta_a F(t, u), \ldots, \Delta_a^k F(t, u), t, u),
\end{equation}
where we have:
\textit{(i)}~$k>0$ and $f, Q$ polynomials,
\textit{(ii)} for~$a\in\Q$, we denote~$\Delta_a:F\in\Q[u][[t]]\mapsto
\frac{F(t, u) - F(t, a)}{u-a}\in\Q[u][[t]]$ and~$\Delta_a^\ell$ the~$\ell^{th}$ iteration
of the \emph{divided difference operator}~$\Delta_a$.
Equations of the form~\eqref{DDE} are called \emph{discrete differential equations (DDEs)}
of~\emph{order}~$k$. 
For instance, \eqref{3const} is a DDE of order~$k=2$, and the divided difference
operator~$\Delta_a$ is taken at the point~$a=1$.
Observe that due to its fixed-point nature, Eq.~\eqref{DDE} admits a unique solution in~$\Q[u][[t]]$
to~\eqref{DDE}. Moreover, the following result proved by Bousquet-Mélou and Jehanne
in~\cite{BMJ06}, reminiscent of Popescu's theorem~\cite[Thm.~$1.4$]{Popescu86},
implies that this solution is annihilated by a nonzero polynomial with coefficients
in~$\Q(t, u)$.
\begin{thm} \emph{(\cite[Thm.~3]{BMJ06})}\label{thm:BMJ06}
Let~$\K$ be a field of characteristic~$0$
and consider two polynomials $f \in \K[u]$ and 
$Q \in \K[x, y_1, \ldots, y_k, t, u]$,
where $k\in\mathbb{N}\setminus \{ 0 \}$.
Let $a\in\K$ and
$\Delta_a :\K[u][[t]]\rightarrow \K[u][[t]]$ be the
divided difference operator $\Delta_a F(t, u) := (F(t, u) - F(t, a))/(u-a)$.
Let us denote by $\Delta_a^{\ell}$ the operator obtained by iterating $\ell$ times $\Delta_a$.
Then, there exists a unique solution $F\in \K[u][[t]]$ to the functional equation Eq.~\eqref{DDE},
and moreover $F(t,u)$ is algebraic over $\K(t, u)$.
\end{thm}
Applying Theorem~\ref{thm:BMJ06} to our example of maps enumeration, there exists
some polynomial~$R\in\Q[t, z]\setminus\{0\}$ such that~$R(t, G(t))=0$.
Namely, $G(t) = 1+t+6t^2+54t^3+594t^4+\cdots$ is annihilated
by~$R(t, z) = 81t^2z^3-9t(9t-2)z^2+(27t^2-66t+1)z-3t^2+47t-1$.
Moreover~$G$ is the unique root of~$R$ in~$\Q[[t]]$.
Thus~$R$ provides a finite amount of data that encodes the sequence~$(c_n)_{n\in\mathbb{N}}$.
Now by writing the linear differential equation satisfied by~$F(t, 1)$ (using~$R$) and by solving the associated recurrence relation, it is possible to show that~$c_0 = 1$
and that~$c_n = \frac{4\cdot 3^{n-1}}{(2n+2)(2n+1)}\binom{3n}{n}$, for all~$n\geq 1$. This is done by 
Additionally one proves,
using the polynomial~$R$, that for all~$n\in\mathbb{N}\setminus\{0, 1\}$
\[(81n^2+81n+18)\cdot c_n-(4n^2+14n+12)\cdot c_{n+1}=0,\]
with~$c_0 = 1, c_1 = 1$. Using the above explicit formula,
computing~$c_{3\cdot 10^6}$ only takes~$5$ seconds in Maple.
Using the above closed formula, one deduces the asymptotic~behavior
\[
c_n\;\sim\;
3^{4n}2^{-2n}(4\sqrt{3\pi}n^{-\frac{5}{2}}+O(n^{-\frac{7}{2}})), \text{ when } n\rightarrow\infty.\]

 \indent
 The oracle that gave us the polynomial~$R(t, z) = 81t^2z^3-9t(9t-2)z^2+(27t^2-66t+1)z-3t^2+47t-1$
 hence allowed us to answer all the questions
 stated at the beginning of this section. We shall emphasize that it is not
 always possible to deduce closed-formulas for such sequences. However, terms of sequences
 satisfying linear recurrences can always be computed very efficiently since the generating functions
 are algebraic, thus solution of a linear
 differential equation with polynomial coefficients~\cite[Prop.~$6.4.3$, Thm.~$6.4.6$]{Stanley99}.\\
 
 A central question is 
 thus to \emph{solve}
 equations like~\eqref{DDE} that is, to compute annihilating polynomials
 for their solutions.
 In various extensions of Theorem~\ref{thm:BMJ06}, the problem considered in this section is in
 spirit the following:
 \begin{center}
 \textit{ For a solution~$F(t, u)$ of a DDE of the form~\eqref{DDE},
     compute~$R\in\K[t, z]\setminus\{0\}$ such that~$R(t, F(t, a))=0$.}\\
   \end{center}

\paragraph{\textbf{Previous algorithmic works. }}
There exists a rich literature regarding the
 effective resolution
 of DDEs. The articles~\cite{BMJ06,BoChNoSa22,BoNoSa23} already contain
 a complete state-of-the-art, but we shall however give a quick overview of this literature.
 Let us start with DDEs such that the polynomial~$Q$ in~\eqref{DDE}
 has degree~$2$ in~$x$. When~$k=1$, Brown introduced in~\cite{Brown65} what is now
 called the~\emph{quadratic method}. This method was generalized later
 by Bender and Canfield in~\cite{BaCa94} to a particular
 family of DDEs of arbitrary order, still with~$Q$ of degree~$2$ in~$x$.\\

 Also, the case where~$Q$ is linear in~$x, y_1, \ldots, y_k$ appears in many walks enumeration
 problems (e.g.~\cite[\S$3.1$]{BMJ06}).
 Solving~DDEs in this case is usually done by applying the~\emph{kernel method},
 introduced by Knuth in~\cite{K68} and further studied by Banderier and Flajolet in~\cite{BaFl02}.
 Since the work of Bousquet-Mélou and Petkov\v{s}ek~\cite{BMP00},
 the linear case is considered as understood. \\

 Regarding strategies based on the~\emph{guess-and-prove} paradigm (popularized in~\cite{Polya48}),
 they have been notably studied
 by Zeilberger~\cite{Zeilberger92} in the case~$k=1$,
 and improved in this same case with Gessel in~\cite{GZ14}.\\
 
 Bousquet-Mélou and Jehanne designed in~\cite{BMJ06} a general method which can be seen
 as a generalization of both the kernel and the quadratic methods. Their strategy
 consists in translating the resolution~of~a~DDE into the resolution of a polynomial
 system admitting a solution with~$F(t, a)$ as one of its coordinates.
We refer to their algorithm as
``duplication''.\\

Recently, intensive studies have been undertaken in~\cite{BoChNoSa22} (resp. in~\cite{BoNoSa23}) in the
direction of effectivity in order to design efficient algorithms for solving DDEs of order $k =1$
(resp. of any order). On the one hand, Bostan, Chyzak, Safey El Din and the author
designed two efficient algorithms based on effective algebraic geometry (we refer to these
algorithms as “elimination” and ”geometry”) and on the other hand they used these
first geometry-driven algorithms in a new hybrid guess-and-prove approach (referred
as ``hybrid'').\\

\begin{center}
    \textit{    The algorithms are now ripe to be delivered turnkey to their natural recipients, combinatorialists, via a Maple package and this article.
    }
\end{center}

\paragraph{\textbf{Contributions.}}
The present paper contains three contributions.
First, it describes the
Maple package~\ddesolver\!, dedicated
to solving DDEs of any order that satisfy
the (harmless) assumptions explicited in~\cite{BoNoSa23}.
More explicitly, \ddesolver contains the function \ann\!.
This function takes as input a DDE of the form~\eqref{DDE}
and its order~$k$, and outputs
a nonzero polynomial annihilating the series~$F(t, a)$. Also, some options can be specified,
namely: the \emph{algorithm} that shall be used among~``duplication'', ``elimination'', ``geometry''
and ``hybrid''; and an option~\emph{variable} explained in~Section~\ref{sec:start}.\\

The second contribution is a careful practical comparison of these algorithms resulting in a table
of timings.\\

The third contribution is the resolution 
with~\ann of
the~DDE associated to the enumeration
of~$3$-greedy Tamari intervals (resp. of~$5$-constellations) in~$1$ minute (resp. in~$3$ hours).
It is the first time that these two equations are solved with an
automatized approach. We however underline that these DDEs are solved without
using effective elimination theory
in~\cite{BMC23} for greedy Tamari intervals and in~\cite{BMS00} for constellations.\\

\paragraph{\textbf{Structure of the paper.}} In Section~\ref{sec:modelling}, we recall
why solving a DDE is reduced~\cite[Sec.~$2$]{BMJ06}
to eliminating variables in a system of polynomial equations.
Also, we introduce in Section~\ref{sec:modelling}
different modeling of the underlying geometric problem. These modelings are useful
later in the paper.
In Section~\ref{sec:preliminaries}, we state some preliminaries on effective algebraic geometry
that are used in Section~\ref{sec:explaining_algorithms} when giving an overview on how
the algorithms ``duplication'', ``elimination'', ``geometry''
and~``hybrid'' (resp.~\cite[\S$3$, \S$5$, \S$6$, \S$4$]{BoNoSa23})
work. In Section~\ref{sec:start}, we explain how to get started with~\ddesolver,
and how to use the options of~\ann\!.
Finally in Section~\ref{sec:examples}, we provide some timings that
illustrate the efficiency of~\ann on various DDEs from the
literature~\cite{BMJ06,Bernardi08,BmFuPr11,BMC23}.\\

\paragraph{\textbf{Notations.}} We gather the notations used in this paper.
First, we always denote by~$\K$ a field of characteristic~$0$, by~$\overline{\K}$
its algebraic closure, by~$\Q$ the field of
rational numbers and by~$\mathbb{F}_p$ the finite field with~$p$ elements.
For a fixed positive integer~$k$, we denote by~$\underline{x}$
(resp.~$\underline{u}$ and~$\underline{z}$) the variables~$x_1, \ldots, x_k$ (resp. $u_1, \ldots, u_k$
and~$z_0, \ldots, z_{k-1}$). 
We denote
by~$\K[x_1, \ldots, x_k]$ the ring of polynomials in the variables~$x_1, \ldots, x_k$
with coefficients in~$\K$.
For~$P\in\K[\underline{x}]$, we denote
by~$V(P)\subset\overline{\K}^r$ its zero set. For any
ideal~$\mathcal{I}\subset\K[\underline{x}]$ and any set of
polynomials~$\mathcal{S}\subset\K[\underline{x}]$, we denote
by~$V(\mathcal{I})\subset\overline{\K}^k$ (resp. by~$V(\mathcal{S})\subset\overline{\K}^k$)
the zero set of~$\mathcal{I}$ (resp. of~$\mathcal{S}$).
Still for~$P\in\K[\underline{x}]$, we denote by~$\opname{disc}_{x_j}(P)$ the
discriminant~\cite[$16$., \S$6$, Ch.$3$]{CoLiOSh15}
of~$P$ with respect to~$x_j$ and by~$\operatorname{LC}_{x_j}(P)$ the leading coefficient of~$P$
with respect to~$x_j$.
We denote by~$\overline{\K}[[t^{\frac{1}\star}]]$ the
ring~$\bigcup_{d\geq 1}\overline{\K}[[t^{\frac{1}{d}}]]$ of Puiseux series in~$t$
with positive fractional exponents.

\section{From DDEs to polynomial systems}\label{sec:modelling}

Starting from a DDE of the form~\eqref{DDE}, we
multiply it by the least power of~$(u-a)$
such that we obtain a polynomial functional equation of the form
\begin{equation}\label{eqn:initpoleqn}
  P(F(t, u), F(t, a), \ldots, \partial_u^{k-1}F(t, a), t, u)=0,
\end{equation}
for some polynomial~$P\in\Q[x, z_0, \ldots, z_{k-1}, t, u]$.
Taking the derivative 
of~\eqref{eqn:initpoleqn} with respect to~$u$ yields by the chain rule
\begin{align}\label{eqn:inipoleqn}
  \partial_uF(t, u)\;\;&\;\cdot\;\;
  \partial_xP(F(t, u), F(t, a), \ldots, \partial_u^{k-1}F(t, a), t, u)\\
 \nonumber &+\;  \partial_uP(F(t, u), F(t, a), \ldots, \partial_u^{k-1}F(t, a), t, u)\; =\; 0.
\end{align}
Consider the equation in~$u$ given by
\begin{equation}\label{eqn:derivative}
  \partial_xP(F(t, u), F(t, a), \ldots, \partial_u^{k-1}F(t, a), t, u)=0.
  \end{equation}
Observe that the solutions
of~\eqref{eqn:derivative} that belong to~$\fracpuiseux$ are also 
solutions, by using~\eqref{eqn:inipoleqn}, of the
equation in~$u$ given by~\[\partial_uP(F(t, u), F(t, a), \ldots, \partial_u^{k-1}F(t, a), t, u)=0.\]
Thus any non constant solution in~$\fracpuiseux$ of~\eqref{eqn:derivative}
is a solution in~$u$ of the system of constraints
\begin{equation}\label{eqn:inisyseqn}
  \begin{cases}
  \;\;\;\; P(F(t, u), F(t, a), \ldots, \partial_u^{k-1}F(t, a), t, u)\;&=0,\\
  \partial_xP(F(t, u), F(t, a), \ldots, \partial_u^{k-1}F(t, a), t, u) \;&=0,
  \;\;\;\;\; u(u- a)\neq 0,\\
  \partial_uP(F(t, u), F(t, a), \ldots, \partial_u^{k-1}F(t, a), t, u) \;&= 0.
  \end{cases}
\end{equation}
 We introduce below hypothesis which is necessary to assume in order to make the
ongoing general strategy work. This assumption holds for the examples considered
in~this~paper.
\begin{align}\label{GA}\tag{\textbf{GA}}\!\underline{\text{Global assumption:}}
  &\;\;\text{There exist $k$ distinct
    solutions~$U_1(t), \ldots, U_k(t)\in\overline{\Q}[[t^{\frac{1}{\star}}]]$}\;\\
 \nonumber &\;\;\,\text{to the system of
equations~\eqref{eqn:inisyseqn}. Also, $U_1(t), \ldots, U_k(t)\notin\overline{\Q}$.} 
\end{align}

In~\cite{BoNoSa23}, the authors studied three geometric interpretations of~\eqref{GA} that
we~recall~below. The first one is introduced in~\cite{BMJ06}, while the last two were
introduced respectively in~\cite[\S$5$, \S$6$]{BoNoSa23}. Also, it is a nontrivial  consequence of the
proof of~\cite[Thm.~$3$]{BMJ06} that the
series~$\{\partial_u^{i-1}F(t, a), U_i(t), F(t, U_i(t))\}_{1\leq i \leq k}$
considered in this paper are
elements of~$\algebraicclosure$.\\

\paragraph{\textbf{Duplication approach.}}
This approach was introduced by Bousquet-Mélou and Jehanne in~\cite[Sec.~$2$]{BMJ06} and works
as follows.
It results from~\eqref{GA} that the following relations hold:
\begin{equation*}
 \forall \;1\leq i \leq k, \begin{cases}
  \;\;\;\; P(F(t, U_i(t)), F(t, a), \ldots, \partial_u^{k-1}F(t, a), t, U_i(t))\;&=0,\\
  \partial_xP(F(t, U_i(t)), F(t, a), \ldots, \partial_u^{k-1}F(t, a), t, U_i(t))\;&=0,\\
  \partial_uP(F(t, U_i(t)), F(t, a), \ldots, \partial_u^{k-1}F(t, a), t, U_i(t))\;&= 0,
  \end{cases}
\end{equation*} and
\begin{equation*}
 \prod\limits_{1\leq i< j \leq k}(U_i(t)-U_j(t))\cdot
  \prod\limits_{1\leq i \leq k}U_i(t)\cdot(U_i(t)-a)\neq 0.
  \end{equation*}
In terms of polynomial equations, the above relations are equivalent to saying that
the~\emph{duplicated} polynomial system~$\Sdup$ defined by
{\footnotesize\begin{equation*}
 \forall \;1\leq i \leq k, \begin{cases}
  \;\;\;\; P(x_i, z_0, \ldots, z_{k-1}, t, u_i)\;&=0,\\
  \partial_xP(x_i, z_0, \ldots, z_{k-1}, t, u_i)\;&=0,\\
  \partial_uP(x_i, z_0, \ldots, z_{k-1}, t, u_i)\;&= 0,
 \end{cases}\\
  m\cdot \prod\limits_{1\leq i< j \leq k}(u_i-u_j)\cdot
  \prod\limits_{1\leq i \leq k}u_i\cdot(u_i-a)-1 = 0,
\end{equation*}
}
admits the nontrivial solutions
\begin{equation}\label{eqn:specialization}
  x_i = F(t, U_i(t)), u_i = U_i(t), z_{i-1} = \partial_u^{i-1}F(t, a),
  \text{ for }~1\leq i \leq k.
  \end{equation}
Note that~\eqref{eqn:specialization} uniquely determines the value associated
to the variable~$m$. Also, the polynomial system~$\Sdup$ admits~$3k+1$ equations
and unknowns (t is considered as a parameter); that is, the polynomials in $\Sdup$ are seen as
elements of the ring~$\Q(t)[m, \underline{x}, \underline{z}, \underline{u}]$.
We denote by~$V(\Sdup)$ the solution set of~$\Sdup$
in~$\algebraicclosure^{k+3}$.\\

\paragraph{\textbf{Elimination theory approach.}}
This second approach was introduced in~\cite[Sec.~$5$]{BoNoSa23}.
The idea is that it follows from~\eqref{eqn:inisyseqn} and~\eqref{GA}
that the system of polynomial constraints defined in~$\algebraicclosure[x, u]$
 ($t$ is again considered as a parameter) by:
\begin{equation}\label{eqn:direct_elim_formulation}
 \begin{cases}
  \;\;\;\; P(x, F(t, a), \ldots, \partial_u^{k-1}F(t, a), t, u)\;&=0,\\
  \partial_xP(x, F(t, a), \ldots, \partial_u^{k-1}F(t, a), t, u)\;&=0,
    \;\;\;\;\; u(u- a)\neq 0,\\
  \partial_uP(x, F(t, a), \ldots, \partial_u^{k-1}F(t, a), t, u)\;&= 0,
  \end{cases}
\end{equation}
admits the~$k$ solutions~$(x, u) = (F(t, U_i(t)), U_i(t))\in\algebraicclosure^2$,
for~$1\leq i \leq k$. Note that
these solutions have distinct~$u$-coordinates.
This observation can
be reformulated geometrically in the following way.
Consider the geometric
projection~$\pi:(x, u, \underline{z})\in\algebraicclosure^{k+2}\mapsto
(\underline{z})\in\algebraicclosure^{k}$,
and denote by~$\mathcal{X}\subset\algebraicclosure^{k+2}$ the solution set
of the polynomial constraints ($t$ is present but still considered as a parameter)
\begin{equation}\label{eqn:initsystem}
 \begin{cases}
  \;\;\;\; P(x, z_0, \ldots, z_{k-1}, t, u)\;&=0,\\
  \partial_xP(x, z_0, \ldots, z_{k-1}, t, u)\;&=0,
    \;\;\;\;\; u(u- a)\neq 0.\\
  \partial_uP(x, z_0, \ldots, z_{k-1}, t, u)\;&= 0,
  \end{cases}
\end{equation}
\indent For~$\bm{\alpha}\in\algebraicclosure^k$,
we denote by~$\#_u(\mathcal{X}, \bm{\alpha})$ the number of~$u$-coordinates
in~$\pi^{-1}(\bm{\alpha})\cap\mathcal{X}$ that are not constants and we define
$ \mathcal{F}_k(u, \mathcal{X})  := \{ \bm{\alpha}\in\algebraicclosure^k \;\;|\;\;
  \#_u(\mathcal{X}, \bm{\alpha})\geq k\}$.
Under assumption~\eqref{GA}, it holds that~$(F(t, a), \ldots, \partial_u^{k-1}F(t, a))\in
  \mathcal{F}_k(u, \mathcal{X})$.\\

  \paragraph{\textbf{Geometric approach.}}
This last approach was introduced in~\cite[Sec.~$6$]{BoNoSa23}.
The idea is that in addition to having~\eqref{eqn:inisyseqn}, it follows from~\eqref{GA}
that the system of polynomial constraints defined in~$\algebraicclosure[x, z_1, u]$ 
($t$ is considered as a parameter) by:
\begin{equation}\label{eqn:direct_elim_formulation}
 \begin{cases}
  \;\;\;\; P(x, F(t, a), z_1, \partial_u^2F(t, a), \ldots, \partial_u^{k-1}F(t, a), t, u)\;&=0,\\
  \partial_xP(x, F(t, a), z_1, \partial_u^2F(t, a), \ldots, \partial_u^{k-1}F(t, a), t, u)\;&=0,
  \;\;\;\;\; u(u- a)\neq 0,\\
  \partial_uP(x, F(t, a), z_1, \partial_u^2F(t, a), \ldots, \partial_u^{k-1}F(t, a), t, u)\;&= 0,
  \end{cases}
\end{equation}
admits the~$k$ solutions~$(x, z_1, u) = (F(t, U_i(t)), \partial_uF(t, a),
U_i(t))\in\algebraicclosure^{k+1}$, for~$1\leq i \leq k$.
Note that these solutions are distinct. This observation can
be reformulated in the following way.\\

Denote by~$\check{z_1}$ the set of variables~$z_0, z_2, \ldots, z_{k-1}$.
Consider the geometric projection~$\pi_{\check{z_1}}:(x, u, \underline{z})\in\algebraicclosure^{k+2}
\mapsto (\check{z_1})\in\algebraicclosure^{k-1}$ and define
the set~
{\footnotesize\[\mathcal{S}_k(\mathcal{X}) := \{\bm{\alpha} = (\alpha_0, \ldots, \alpha_{k-1})
\in\algebraicclosure^k \;\,|\,\;
\bm{\alpha}\in\pi(\mathcal{X}) \;\,\wedge\;\, \mathcal{X} \cap\;
\pi_{\check{z_1}}^{-1}((\alpha_0, \alpha_2, \ldots, \alpha_{k-1}))\geq
k\}.\]
}
Under assumption~\eqref{GA}, it holds that~$(F(t, a), \ldots, \partial_u^{k-1}F(t, a))\in
  \mathcal{S}_k(\mathcal{X})$.\\

 \paragraph{\textbf{Summary.}} In Section~\ref{sec:modelling}, we defined for each of the three approaches
  a solution set which contains a point
  whose~$z_0$-coordinate is the series~$F(t, a)$. The goal of~Section~\ref{sec:explaining_algorithms} is
  to explain how the algorithms studied in~\cite[\S$3$, \S$5$, \S$6$]{BoNoSa23}
  compute a polynomial characterization of these sets and deduce from this characterization
  an annihilating polynomial of~$F(t, a)$.

  \section{Some preliminaries on polynomial tools}\label{sec:preliminaries}

  The current section is devoted to the introduction of the polynomial tools used
  in~Section~\ref{sec:explaining_algorithms}.

\subsection{Elimination theory through Gr\"obner bases}\label{sec:GB}

For~$x_1, \ldots, x_n$ some variables and~$\alpha = (\alpha_1, \ldots, \alpha_n)\in\mathbb{Z}_{\geq 0}^n$,
we denote~$\bm{x}^{\bm{\alpha}} := x_1^{\alpha_1}\cdots x_n^{\alpha_n}$.
Also for~$\mathcal{S}\subset\K[x_1, \ldots, x_n]$, we denote by~$\langle \mathcal{S} \rangle$
the ideal generated by~$\mathcal{S}$ in~$\mathbb{K}[x_1, \ldots, x_n]$.
\paragraph{Monomial orders and Gr\"obner bases}
 The idea of elimination monomial orderings is to attribute to a
variable (or to a block of variables) that we want to eliminate
a larger weight than the weights of the other variables. An example is the
 lexicographic order~$\succ_{\text{lex}}$.

\begin{definition}(\cite[Def.$3$, \S$2$, Ch.$2$]{CoLiOSh15})
   Let~$\alpha =(\alpha_1, \ldots, \alpha_n)$ and~$\beta = (\beta_1, \ldots, \beta_n)$ be
  in~$\mathbb{Z}_{\geq 0}^n$. We say that~$\bm{x}^{\bm{\alpha}}\succ_{\text{lex}} \bm{x}^{\bm{\beta}}$
  if the leftmost nonzero entry of the vector difference~$\alpha-\beta\in\mathbb{Z}^n$ is positive.
\end{definition}

\begin{tcolorbox}[colback=green!5!white,colframe=green!75!green,]
\begin{example}
In~$\K[x_1, x_2]$: $x_1^4x_2^2\succ_{\text{lex}}x_1^3x_2^{10}$, $x_1^2x_2\succ_{\text{lex}}x_1$ and~$x_1^2\succ_{\text{lex}} x_2$.
\end{example}
\end{tcolorbox}

Another useful family of monomial orders that we will refer to in~Section~\ref{sec:explaining_algorithms}
is the family of graded monomial orders. An example is the
graded reverse lexicographic order~$\succ_{\text{grevlex}}$.

\begin{definition}(\cite[Def.$6$, \S$2$, Ch.$2$]{CoLiOSh15})
  Let~$\alpha=(\alpha_1, \ldots, \alpha_n)$ and~$\beta=(\beta_1, \ldots, \beta_n)$ be in
  $\mathbb{Z}_{\geq 0}^n$. We say
  that~$\bm{x}^{\bm{\alpha}}\succ_{\text{grevlex}}\bm{x}^{\bm{\beta}}$
  if~$|\alpha| = \sum_{i=1}^n\alpha_i>|\beta| = \sum_{i=1}^n\beta_i$,
  or if~$|\alpha| = |\beta|$ and the rightmost nonzero entry
  of~$\alpha - \beta\in\mathbb{Z}^n$ is negative. 
  \end{definition}

\begin{tcolorbox}[colback=green!5!white,colframe=green!75!green,]
\begin{example}
  In~$\K[x_1, x_2, x_3]$:
  $x_1^5x_2^7x_3\succ_{\text{grevlex}} x_1^4x_2^2x_3^3$, $x_1x_2^4x_3^2\succ_{\text{grevlex}}x_1^3x_2x_3^3$ 
  \end{example}
\end{tcolorbox}

We invite readers
unfamiliar with monomial orders (see~\cite[Def.$1$,\S$2$,Ch.$2$]{CoLiOSh15} for a general definition)
to choose one of the two above monomial orders~$\succ_{\text{lex}}$
or~$\succ_{\text{grevlex}}$ for the next definition.
Also, for~$p\in\mathbb{K}[x_1, \ldots, x_n]$ a polynomial and~$\succ$ a monomial order
on~$\K[x_1, \ldots, x_n]$, we denote by~$\operatorname{LT}_{\succ}(p)$ the leading term of~$p$
with respect to the monomial order~$\succ$.

\begin{definition}(\cite[Def.$5$, \S$5$, Ch.$2$]{CoLiOSh15})\label{def:gb}
  Fix a monomial order~$\succ$ on the polynomial ring~$\K[x_1, \ldots, x_n]$. A finite
  subset~$G=\{g_1, \ldots, g_s\}$ of an ideal~$\mathcal{I}\subset\K[x_1, \ldots, x_n]$ different
  from~$\{0\}$ is said to be a Gr\"obner basis of~$\mathcal{I}$ for the order~$\succ$ if
 $\langle \LT(g_1), \ldots, \LT(g_n) \rangle = \langle \{\LT(g) \,|\,
  g\in \mathcal{I}\} \rangle$.
\end{definition}

A fundamental property~\cite[Cor.$6$, \S$5$, Ch.$2$]{CoLiOSh15}
of Gr\"obner bases is, with the notations of Definition~\ref{def:gb}, that
such a basis~$G$ always exists, and that~$G$ generates~$\mathcal{I}$, that
is~$\langle G \rangle = \mathcal{I}$.\\

  \paragraph{\textbf{Gr\"obner bases and projections.}}

  The results below justify the use of Gr\"obner bases as a dedicated theoretical and computational
  tool, allowing one to characterize projections as solution sets of
  conjunctions of polynomial equations and inequations.
  
  \begin{thm}[Elimination theorem](\cite[Thm.$2$, \S$1$, Ch.$3$]{CoLiOSh15})\label{thm:elimination}
    Let~$\mathcal{I}\subset\mathbb{K}[x_1, \ldots, x_n]$ be an ideal.
    Denote by~$G$ a Gr\"obner basis of~$\mathcal{I}$ with respect to the 
    order~$\succ_{\text{lex}}$. Then for every~$0\leq \ell \leq n-1$, the set~$G_\ell =
    G\cap\mathbb{K}[x_{\ell+1}, \ldots, x_n]$ is a Gr\"obner basis of the
    ideal~$\mathcal{I}_\ell := \mathcal{I}\cap\mathbb{K}[x_{\ell+1}, \ldots, x_n]$.
  \end{thm}
  By Theorem~\ref{thm:elimination},
  the set of polynomials~$G_\ell$ finitely generates~$\mathcal{I}_\ell$. Let us now introduce
  the geometric projection~$\pi_\ell:(x_1, \ldots, x_n)\in\overline{\K}^n\mapsto(x_{\ell+1}, \ldots, x_n)
  \in\overline{\K}^{n-\ell}$. Recall that a Zariski closed set in~$\overline{\K}^n$ is defined as
  the solution set of some polynomial equations defined in~$\K[x_1, \ldots, x_n]$.

  \begin{thm}[Closure theorem](\cite[Thm.$3$, \S$2$, Ch.$3$]{CoLiOSh15})\label{thm:closure}
    Let~$\mathcal{I}\subset\K[x_1, \ldots, x_n]$ be an ideal and~$V(\mathcal{I})\subset\overline{\K}^n$
    its zero set. Consider~$\mathcal{I}_\ell := \mathcal{I} \cap\K[x_{\ell+1}, \ldots, x_n]$.
    Then~$V(\mathcal{I}_\ell)$ is the smallest Zariski closed set containing~$\pi_\ell(V(\mathcal{I}))
    \subset\overline{\K}^{n-\ell}$.
  \end{thm}
  
  Theorems~\ref{thm:elimination} and~\ref{thm:closure}
  imply that characterizing the Zariski closure of the projection of a zero set~$V(\mathcal{I})$
  onto some coordinate subspace is done by computing a Gr\"obner basis for~$\succ_{\text{lex}}$.
 
  \begin{thm}[Extension theorem](\cite[Thm.$2$, \S$5$, Ch.$3$]{CoLiOSh15})\label{thm:extension}
    Let~$G=\{g_1, \ldots, g_s\}$ be a Gr\"obner basis of~$\mathcal{I}\subset\K[x_1, \ldots, x_n]$
    for the order~$\succ_{\text{lex}}$. For each~$1\leq j\leq s$, consider
    \[ g_j = c_j(x_2, \ldots, x_n)\cdot x_1^{N_j} + \text{(terms in which~$x_1$ has degree~$<N_j$)},\]
    where~$N_j\geq 0$ and~$x_j\in\K[x_2, \ldots, x_n]$ is nonzero.
    Assume~$\bm{\alpha} = (\alpha_2, \ldots, \alpha_n) \in V(\mathcal{I}\cap\K[x_2, \ldots, x_n])$ is
    a partial solution with the property that~$\bm{\alpha}\notin V(c_1, \ldots, c_s)$.
    Then
    \[ \{f(x_1, \bm{\alpha}) \, | \, f\in\mathcal{I} \} =
    \{ g_{j_0}(x_1, \bm{\alpha})\},\]
    where~$g_{j_0}\in G$ satisfies~$c_{j_0}(a)\neq 0$ and~$g_{j_0}$ has minimal degree in~$x_1$
    among all elements~$g_j$ with~$c_j(\bm{\alpha})\neq 0$. Furthermore,
    if~$g_{j_0}(\alpha_1, \bm{\alpha})=0$ for~$\alpha_1\in\overline{\K}$,
      then~$(\alpha_1, \bm{\alpha})\in V(\mathcal{I})$.
    \end{thm}

  A consequence of Theorem~\ref{thm:extension} is that characterizing projections
  (and not their Zariski closure)
  can be done by considering the
  disjunction of inequations associated to the non simultaneous vanishing
  of the~$c_1, \ldots, c_s$ of Theorem~\ref{thm:extension}: for each condition in this disjunction,
  one shall add the vanishing of all the elements of~$G_\ell$ of Theorem~\ref{thm:elimination}.\\

  \begin{tcolorbox}[colback=green!5!white,colframe=green!75!green,]
  \begin{example}
    If~$s=2$, $G_\ell = \{g_1, g_2, g_3\}$,
    the projection is the solution set of the constraints\vspace{-0.2cm}
\[   (c_1 \neq 0 \; \wedge\; f_1=0 \;\wedge\; f_2 = 0 \;\wedge\; f_3 = 0)  \;\vee
\;  (c_2 \neq 0 \;\wedge\; f_1=0 \;\wedge\; f_2 = 0 \;\wedge\; f_3 = 0).\]
    \end{example}
\end{tcolorbox}
  
\subsection{Counting specific solutions}
Let~$n, s$ be positive integers. We take
the notation~$\bm{x} = x_1, \ldots, x_n$ and~$\bm{y} = y_1, \ldots, y_s$.
Let~$\mathcal{I}\subset\K[\bm{x}]$ be a radical ideal such that its zero
set~$V(\mathcal{I})\subset\overline{\K}^n$ is finite.
For~$\bm{\alpha}\in V(\mathcal{I})$, we denote by~$[x_1](\bm{\alpha})$ the~$x_1$-coordinate
of~$\bm{\alpha}$. We introduce the
projection~$\pi_{x_1}:(\bm{x})\in\overline{\K}^n\mapsto x_1\in\overline{\K}$.
 Let~$g\in\K[\bm{y}][z]$ and
denote $\opname{LC}_z(g)\in\K[\bm{y}]$ the leading coefficients of~$g$ in~$z$.\\

 This current subsection introduces tools for answering the two problems below.
 
\begin{prb}\label{qst1}
  Let~$\ell$ be a positive integer. 
  Characterize with polynomial equations defined in~$\K[\bm{x}]$
  the set~$\{\bm{\beta}\in V(\mathcal{I}) \, |\,
  \pi_{x_1}^{-1}([x_1](\bm{\beta}))\cap V(\mathcal{I}) \geq \ell\}$.
\end{prb}
  
\begin{prb}\label{qst2}
  Let~$1\leq \ell \leq \deg_z(g)$. Characterize with polynomial inequations defined in~$\K[\bm{y}]$
  the points~$\bm{\beta}\in\overline{\K}^s\setminus V(\opname{LC}_x(g))$ such
  that~$g(\bm{y} = \bm{\beta}, z)$ has at least~$\ell$ distinct roots.
\end{prb}

\subparagraph{Answer to Problem~\ref{qst1}}
Denote~$\mathcal{A} := \K[\bm{x}]/\mathcal{I}$. It results from the finiteness
of~$V(\mathcal{I})$ that~$\mathcal{A}$ has finite dimension as a~$\K$-vector
space~\cite[Thm.$6$, \S$3$, Ch.$5$]{CoLiOSh15}.
Define~$m_{x_1}:f\in\mathcal{A}\mapsto x_1\cdot f\in\mathcal{A}$
to be the~\emph{multiplication map by~$x_1$} in~$\mathcal{A}$, and denote by~$\chi_{x_1}\in\K[T]$ its
characteristic polynomial.
The following results from the radicality of~$\mathcal{I}$ and
from~\cite[Prop.$2.7$, \S$2$, Ch.$4$]{UAG}.

\begin{fact}\label{factStick}
  We have the equality~$\chi_{x_1} = \prod_{\bm{\alpha}\in V(\mathcal{I})} (T - [x_1](\bm{\alpha}))$.
  \end{fact}
As a consequence of Fact~\ref{factStick}, we obtain the following set equality:
\begin{align}\label{eqn:qst1}
\{\bm{\beta}\in V(\mathcal{I}) &\, |\,
\pi_{x_1}^{-1}([x_1](\bm{\beta}))\cap V(\mathcal{I}) \geq \ell\}
\\
\nonumber &= \{ \bm{\beta} \in \overline{\K}^n \, |\,
  \bm{\beta}\in V(\mathcal{I}) \wedge \chi_{x_1}([x_1](\bm{\beta}))=0 \wedge \cdots
  \wedge \partial_T^{\ell-1}\chi_{x_1}([x_1](\bm{\beta}))=0\}.
\end{align}
Equality~\eqref{eqn:qst1} provides a conjunction of polynomial equations that answers
Problem~\ref{qst1}.\\

\begin{tcolorbox}[colback=green!5!white,colframe=green!75!green,]
\begin{example}
  Assume~$G = \{g_1, g_2\}\subset\K[x_1, x_2]$ generates~$\mathcal{I}$ and~$\ell = 2$.
  The conjunction of polynomial equations in~$\K[x_1, x_2]$ is
  $ g_1(x_1, x_2)=0 \wedge g_2(x_1, x_2)=0 \wedge \chi_{x_1}(x_1) = 0 .$
\end{example}
\end{tcolorbox}\vspace{-0.2cm}

\subparagraph{Answer to Problem~\ref{qst2}}

Define the~\emph{Hermite quadratic form} associated with~${g}$
by~\[H_g:(f, h)\in (\K(\mb{y})[z]/\langle g\rangle )^2\mapsto
\opname{Trace}(m_{f\cdot h})\in\mathbb{K}(\mb{y}),\]
where~$\opname{Trace}(\cdot)$ is the trace operator and~$m_{f\cdot h}$ is the multiplication map
by~$f\cdot h$ in~$\K(\mb{y})[z]/\langle g\rangle $. Also, denote
by~$M_{H_g}\in\K(\bm{y})^{\deg_z(g)\times\deg_z(g)}$
the matrix of~$H_g$ in the basis~$\{1, z, \ldots, z^{\deg_z(g)-1}\}$.\\

The following fact is an immediate consequence
of~\cite[Thm.$5.2$, \S$5$, Ch.$2$]{UAG} and of the observation that,
by perfoming euclidian divisions by~$g$
in~$\K(\bm{y})[z]$, the denominator
of the image of~$z^{i}$ in~$\K(\mb{y})[z]/\langle g\rangle $ can only be a power of~$\opname{LC}_z(g)$.

\begin{fact}\label{factHermite}
  The denominators in the matrix~$M_{H_g}$ are powers
  of~$\opname{LC}_z(g)$. Moreover if~$1\leq \ell \leq \deg_z(g)$,
   then the points
  $\bm{\beta}\in\overline{\K}^s\setminus V(\opname{LC}_z(g))$ at which
   $g(\bm{y} = \bm{\beta}, z)$ admits less than~$\ell$ distinct solutions are
   precisely the points~$\bm{\beta}\in\overline{\K}^s\setminus V(\opname{LC}_z(g))$
 at which the~$\ell\times\ell$-minors of~$M_{H_g}$ all vanish.
  \end{fact}

Fact~\ref{factHermite} implies that the
conjunction of polynomial inequations answering Problem~\ref{qst2} is given by the 
non vanishing of~$\opname{LC}_z(g)$ and 
by the non vanishing of the~$\ell\times\ell$-minors~of~$M_{H_g}$.\\

\begin{tcolorbox}[colback=green!5!white,colframe=green!75!green,]
\begin{example}
  Assume that the~$\ell\times \ell$-minors of~$M_{H_g}$ are~$m_1, m_2, m_3, m_4
  \in\K(\mb{y})$.
  The conjunction of inequations in~$\K(\bm{y})$ is\vspace{-0.2cm}
  \[ \opname{LC}_z(g) \neq 0 \wedge m_1(\bm{y})\neq 0 \wedge m_2(\bm{y})\neq 0
  \wedge m_3(\bm{y}) \neq 0
  \wedge m_4(\bm{y}) \neq 0.\]
\end{example}
\end{tcolorbox}

\subsection{Change of monomial ordering}

For $n$ generic elements of~$\K[x_1, \ldots, x_n]$ of degree $d$,
computing a Gr\"obner basis of~$\mathcal{I}$ for the order~$\succ_{\text{grevlex}}$
has an arithmetic cost which is in~$d^{O(n)}$~\cite[Prop.$1$]{BaFaSa15}.\\

Let~$\mathcal{I}\subset\K[x_1, \ldots, x_n]$ be an ideal whose solution set is finite.
The following fact mentions an algorithm called~FGLM~\cite{FGLM93}
that on input a Gr\"obner basis of~$\mathcal{I}$ for the order~$\succ_{\text{grevlex}}$
outputs a Gr\"obner basis
of~$\mathcal{I}$ for the order~$\succ_{\text{lex}}$.
The complexity of this algorithm is, under some genericity assumptions,
in~$\tilde{O}(D^\omega)$\footnote{We denote
by~$\omega$ the constant of multiplication matrix that is,
multiplying two matrices in~$\K^{n\times n}$ can be done by using~$\tilde{O}(n^\omega)$
arithmetic operations in~$\mathbb{K}$, where~$\tilde{O}(\cdot)$ is a~$O(\cdot)$ that
hides polylogarithmic factors.} (see~\cite{NeSc20}), where~$D$ is the cardinality of~$V(\mathcal{I})$
($\mathcal{I}$ is assumed radical).
By the Bézout bound, and with the notations of the above paragraph, we have the bound~$D\leq d^n$,
this makes the following fact useful for computing in~$d^{O(n)}$ a Gr\"obner bases of~$\mathcal{I}$
for a~$\succ_{\text{lex}}$ order.

\begin{fact}\label{fact:fglm}
  On input an ideal~$\mathcal{I}\subset\K[x_1, \ldots, x_n]$ whose solution
  set~$V(\mathcal{I})\subset\overline{\K}^n$ is finite, there exists an algorithm
  that computes a Gr\"obner basis~$G_{\succ_{\text{grevlex}}}$
  of~$\mathcal{I}$ for the order~$\succ_{\text{grevlex}}$ and that
  turns it into a Gr\"obner basis~$G_{\succ_{\text{lex}}}$ of~$\mathcal{I}$ for the order~$\succ_{\text{lex}}$.
  \end{fact}

We mention that with the notations as above and without using the algorithm FGLM~\cite{FGLM93}
mentioned in~Fact~\ref{fact:fglm}, the computation of a Gr\"obner basis of~$\mathcal{I}$ for
the order~$\succ_{\text{lex}}$ has an arithmetic cost bounded
by $C_2d^{C_3n^3}$~\cite{CaGaHe88}, for~$C_2, C_3\in\mathbb{Z}_{\geq 0}$.

  \section{Function implemented in~\ddesolver}\label{sec:explaining_algorithms}

  The function implemented in~\ddesolver \, is:
  \begin{center}\ann\end{center}

We explain in this section the algorithms on which this function, and its options, rely.\\

 The function~\ann  computes an annihilating polynomial of the series~$F(t, a)$. Four algorithms can be applied:
 ``duplication'', ``elimination'', ``geometry'', ``hybrid''. These algorithms are studied
 respectively in~\cite[\S$3$, \S$5$, \S$6$, \S$4$]{BoNoSa23}.
 As the algorithm~``hybrid'' relies on the
 three other algorithms, we first describe~``duplication'', ``elimination'' and~``geometry''.\\
 
 In~Section~\ref{sec:modelling}, we defined solution sets that contain a point
  whose~$z_0$-coordinate is~$F(t, a)$:
  \begin{itemize}
  \item ``duplication'': $V(\Sdup)$ contains the point defined
    by~\eqref{eqn:specialization},
      \item ``elimination'': $\mathcal{F}_k(u, \mathcal{X})$ contains the
        point~$(F(t, a), \ldots, \partial_u^{k-1}F(t, a))$,
      \item ``geometric'': $\mathcal{S}_k(\mathcal{X})$ contains the
        point~$(F(t, a), \ldots, \partial_u^{k-1}F(t, a))$.
      \end{itemize}
 \indent  We introduce a finiteness assumption that holds for the examples studied in this paper.
  \begin{align}\label{FA}\tag{\textbf{FA}}\underline{\text{Finiteness assumption:}}
    &\;\;\text{The sets
      $V(\Sdup)$, $\mathcal{F}_k(u, \mathcal{X})$ and~$\mathcal{S}_k(\mathcal{X})$
    are finite.}
\end{align}

  \paragraph{\textbf{General strategy.}} Under~\eqref{FA}, the general spirit of the algorithms
  from~\cite[\S$3$, \S$5$, \S$6$]{BoNoSa23} is:
  \begin{enumerate}
\item[(i)] To compute a disjunction of
  polynomial equations whose solution set is the set of interest,
\item[(ii)] To eliminate all variables except~$(z_0, t)$ from the polynomial characterization
  of step~$\text{(i)}$.
    \end{enumerate}

We apply the above steps to the three approaches from~Section~\ref{sec:modelling}.

  \subsection{Option ``duplication''}\label{sec:idea_duplication}

  \;\;\;\;\, \textit{Step~\textbf{(i)} of the general strategy:}
  The polynomial system~$\Sdup$
  characterizes the set~$V(\Sdup)$.\\
  
  \textit{Step~\textbf{(ii)} of the general strategy:}
Formally, we denote by~$\mathcal{I}_{\text{dup}}\subset\Q(t)[m, \underline{x}, \underline{z},
    \underline{u}]$ the ideal generated by the polynomials in~$\Sdup$.
  By assumption~\eqref{FA}, it is possible to apply to the ideal~$\mathcal{I}_{\text{dup}}$
  the algorithm from Fact~\ref{fact:fglm}.
  This algorithm outputs,
  for a proper choice of variable ordering, a Gr\"obner basis~$G_{\succ_{\text{lex}}}$
  of~$\mathcal{I}_{\text{dup}}$ that contains an element
  of~$\Q(t)[z_0]$. Applying Theorem~\ref{thm:closure} to the projection onto the~$z_0$-coordinate
  space together with Theorem~\ref{thm:elimination} yields, under~\eqref{FA}, 
  that this polynomial is nonzero.
  Denote by~$R\in\Q[t, z_0]$ its numerator. 
  By~\cite[Prop.~$2$]{BoNoSa23}, we have~$R(t, F(t, a))~=~0$.\\

  We refer the reader to Annex~\ref{annex:dup} for the resolution of~$3$-constellations via this method.

\subsection{Option ``elimination''}\label{sec:elim}

Recall from Section~\ref{sec:modelling}
that~$\mathcal{F}_k(u, \mathcal{X}) = \{ \bm{\alpha}\in\overline{\Q(t)}^{k} \, | \,
\#_u(\mathcal{X}, \bm{\alpha})\geq k\}$. 
The algorithm from~\cite[\S$5$]{BoNoSa23} requires the additional technical
assumption that for any~$\bm{\alpha}\in\overline{\Q(t)}^k$,
we have~$\#_u(\mathcal{X}, \bm{\alpha})<+\infty$.
We do not make it explicit for
sake of simplicity. Under this technical assumption and~\eqref{FA},
the algorithm designed in~\cite[\S$5$]{BoNoSa23} works as follows.\\

\textit{Step~\textbf{$\mathbf{(i)}$} of the general strategy:}
\begin{enumerate}
\item We compute using Theorems~\ref{thm:elimination} and~\ref{thm:extension}
   successive disjunctions of conjunctions
  of polynomial equations and inequations defined in~$\K(t)[u, \underline{z}]$
  whose solution sets are successively the projections onto the~$(u, \underline{z})$-coordinate space,
  then onto the~$(\underline{z})$-coordinate space. \\
  The union of the polynomial constraints in these disjunctions have the form $\{I, E\}$,
  where~$I$ is a set of inequations defined in~$\Q(t)[u, \underline{z}]$
  and~$E$ is a set of equations defined in~$\Q(t)[u, \underline{z}]$.
  Without loss of generality, we assume that the
  polynomials in~$E$ are a Gr\"obner basis of the ideal they generate,
  for the order~$\succ_{\text{lex}}$ with~$u$ greater than~$\underline{z}$.
\item We apply~Theorem~\ref{thm:extension} to~$E$:
  at~$\underline{z} = \bm{\alpha}\in\overline{\Q(t)}^k$ fixed,
  the cardinality condition in the definition of~$\mathcal{F}_k(u, \mathcal{X})$ is equivalent to
  studying for a polynomial~$g_{j_0}$ in~$E$:
  \begin{itemize}
  \item if $g_{j_0}$ has its leading coefficient in~$u$ that does not vanish
    at~$\underline{z} = \bm{\alpha}$ and, by~Theorem~\ref{thm:extension},
    if~$g_{j_0}$ is the polynomial in~$E$ of minimal degree in~$u$
    that satisfies this property: checking this minimality condition
    yields a disjunction of polynomial
    equations and inequations defined in~$\Q(t)[\underline{z}]$,
  \item for such a~$g_{j_0}$, it remains to add a conjunction of polynomial conditions
    in~$\Q(t)[u, \underline{z}]$ so
    that~$g_{j_0}(u, \underline{z} = \bm{\alpha})$ has at least~$k$ distinct solutions: \\
    If~$\deg_{u}(g_{j_0})<k$, adds the vanishing of all the coefficients (in~$u$) of~$g_{j_0}$ to~$E$;
    else adds the conjunction of inequations given by the non vanishing of the~$k\times k$-minors
    of~$M_{H_{g_{j_0}}}$ from Fact~\ref{factHermite}.
  \end{itemize}
\item Steps~$\mathbf{1}$ and~$\mathbf{2}$ above
  computed a disjunction of polynomial equations and inequations
  defined in~$\Q(t)[u, \underline{z}]$ whose solution set in~$\overline{\Q(t)}^{k+1}$
  is~$\mathcal{F}_k(u, \mathcal{X})$. As by~\eqref{FA} this set is finite,
  this polynomial characterization can be turned into a disjunction of conjunction of polynomial
  equations by using Rabinowitsch trick\footnote{Instead of considering an inequation~$Q\neq 0$ (for a
  given polynomial~$Q\in\K[x_1, \ldots, x_n]$), we introduce an extra variable~$m$ and
  consider the equation~$m\cdot Q-1=0$ defined in~$\K[m, x_1, \ldots, x_n]$.}
  in order to remove the solution set given by the inequations:
  the introduced variables are eliminated using~Theorem~\ref{thm:elimination}. We denote
  by~$\mathcal{D}_{\text{elim}}$ the so-obtained disjunction.
  \end{enumerate}

\textit{Step~\textbf{(ii)} of the general strategy:}
We perform this step in the same spirit as what we did in~Section~\ref{sec:idea_duplication}.
By~\eqref{FA}, the solution set of~$\mathcal{D}_{\text{elim}}$ is finite. We can thus apply
the algorithm from Fact~\ref{fact:fglm} to each disjunction in~$\mathcal{D}_{\text{elim}}$.
Each application of this algorithm allows to compute a nonzero polynomial
of~$\Q(t)[z_0]$: this polynomial might be equal to~$1$ when there is no solution to the studied
conjunction. We consider the numerators of all these polynomials and denote by~$R\in\Q[t, z_0]$
their product. By~\cite[Prop.~$5.5$]{BoNoSa23}, we have~$R(t, F(t, a))~=~0$.\\

  We refer the reader to Annex~\ref{annex:elim} for the resolution of~$3$-constellations via this method.

\subsection{Option ``geometry''}\label{sec:geometry}

Recall from Section~\ref{sec:modelling}
that~
{\footnotesize\[\mathcal{S}_k(\mathcal{X}) := \{\bm{\alpha} = (\alpha_0, \ldots, \alpha_{k-1})
\in\algebraicclosure^k \;\,|\,\;
\bm{\alpha}\in\pi(\mathcal{X}) \;\,\wedge\;\, \mathcal{X} \cap\;
\pi_{\check{z_1}}^{-1}((\alpha_0, \alpha_2, \ldots, \alpha_{k-1}))\geq
k\}.\]
}
The algorithm from~\cite[\S$6$]{BoNoSa23} requires many additional technical
assumptions~(see the beginning of~\cite[\S$6$]{BoNoSa23}).
We do not make them explicit for
sake of simplicity.
In the rest of this section, denote by~$\succ_{\text{bgrevlex}}$ the block monomial order
over~$\Q(t)[m, x, u, \underline{z}]$ defined as follows:
\begin{itemize}
\item We use the monomial order~$\succ_{\text{grevlex}}$ on each of the two
  blocks~$\{m, x, u, z_1\}$ and~$\{\check{z_1}\}$,
\item Two monomials are compared w.r.t. the variables~$m, x, u, z_1$. In case of equality,
  they are compared w.r.t. the variables~$\check{z_1}$.
\end{itemize}

\begin{tcolorbox}[colback=green!5!white,colframe=green!75!green,]
\begin{example}
  For~$k=3$:\vspace{-0.2cm}
  \[x^3uz_1^2m \succ_{\text{bgrevlex}} z_0^{20}z_2^3,\;\;
  x^3uz_1^2mz_2 \succ_{\text{bgrevlex}} x^3uz_1^2mz_0 \;\;\text{  and }~z_0^4\succ_{\text{bgrevlex}}z_0^2z_2.\]
  \end{example}
\end{tcolorbox}\vspace{0.2cm}

Under these assumptions and~\eqref{FA},
the algorithm designed in~\cite[\S$6$]{BoNoSa23} works as follows.\\

\textit{Step~\textbf{$\mathbf{(i)}$} of the general strategy:}

\begin{enumerate}
\item Compute a Gr\"obner basis~$G$ of the ideal~$\mathcal{I} := \langle P, \partial_xP, \partial_uP,
  m\cdot u(u-a)-1\rangle\subset\Q(t)[m, x, u, \underline{z}]$
  for the order~$\succ_{\text{bgrevlex}}$.
\item Compute, using normal form computations
  modulo~$\mathcal{I}$\footnote{A normal form computation modulo
an ideal is the multivariate generalization
of an euclidean division by a polynomial.}, the matrix~$M_{z_1}$ of the multiplication map
  $m_{z_1}$ in the quotient ring~$\Q(t, \check{z_1})[m, x, u, z_1]/j(\mathcal{I})$,
  where $j(\mathcal{I})$ is the image of the injective map~$j:\Q(t)[m, x, u, \underline{z}]
  \rightarrow\Q(t, \check{z_1})[m, x, u, z_1]$.
\item Compute the characteristic polynomial~$\chi_{z_1}\in\Q(t, \check{z_1})[T]$.
  Using the assumptions of~\cite[\S$6$]{BoNoSa23}, we can replace~$\chi_{z_1}$ by its numerator
  with respect to the variables~$\check{z_1}$, so that we now have~$\chi_{z_1}\in\Q(t)[\check{z_1}, T]$.
\item We denote by~$\mathcal{D}_{\text{geom}}$ the set of equations given by the vanishing of the
  polynomials in~$G$ and of the new
  conditions~$(\chi_{z_1})|_{T=z_1}=0 \wedge \ldots \wedge  (\partial_{T}^{k-1}\chi_{z_1})|_{T=z_1}=0$.
  \end{enumerate}

\textit{Step~\textbf{$\mathbf{(ii)}$} of the general strategy:}
We perform this step in the same spirit as we did in~Section~\ref{sec:idea_duplication}.
By~\eqref{FA}, the solution set of~$\mathcal{D}_{\text{geom}}$ is finite. We thus apply
the algorithm from Fact~\ref{fact:fglm} to the ideal generated by the polynomials associated with
the equations in~$\mathcal{D}_{\text{geom}}$.
The application of this algorithm allows to compute a nonzero polynomial
of~$\Q(t)[z_0]$. We consider the numerator of this polynomial and denote it by~$R\in\Q[t, z_0]$.
By~\cite[Prop.~$6.4$]{BoNoSa23}, we have~$R(t, F(t, a))=0$.\\

  We refer the reader to Annex~\ref{annex:geom} for the resolution of~$3$-constellations via this method.

\subsection{Option ``hybrid''}\label{sec:hybrid}

The \emph{hybrid guess-and-prove} algorithm~\cite[\S$4$]{BoNoSa23} works as follows:
\begin{enumerate}
\item  Compute bounds~$b_{z_0}, b_t$ such that~$\deg_t(R)\leq b_t$
  and~$\deg_{z_0}(R)\leq b_{z_0}$, for some nonzero polynomial~$R\in\Q[t, z_0]$ annihilating~$F(t, a)$
  (we detail this step in the next section),
\item Compute the truncated series~$F(t, a)\bmod t^{2b_tb_{z_0}+1}$,
\item Guesses a polynomial~$M\in\Q[t, z_0]$ such that~$M(t, F(t, a)) = O(t^{(b_t+1)(b_{z_0}+1)-1})$,
  \item Checks that~$M(t, F(t, a))=O(t^{2b_tb_{z_0}+1})$.
  \end{enumerate}

When step~$\mathbf{4}$ is satisfied, applying~\cite[Prop.$5$]{BoNoSa23}
yields~$M(t, F(t, a))=0$.\\

We refer the reader to Annex~\ref{annex:hybrid} for the resolution of~$3$-constellations
via this method.

\subsection{Details of implementations}

As we target practical efficiency,
the implementations of the options ``duplication'', ``elimination'' and~``geometry''
incorporate the following key improvements from computer algebra. Recall that
the output of these three methods is a nonzero~$R\in\Q[t, z_0]$ such that~$R(t, F(t, a))=0$.\\

To reduce the computations of~$\Q(t)$~to~$\Q$, we perform~evaluation-interpolation on~$t$.\\

Also, we rely on fast multi-modular arithmetic. This consists in applying multiple times
a given algorithm
with the base field~$\mathbb{Q}$ replaced by successive distinct prime fields~$\mathbb{F}_p$.
These computations output
 various images mod~$p$ of the same~$R\in\Q[t, z_0]$.
From there, one lifts the modular coefficients over~$\Q$ by applying
the Chinese Remainder Theorem (CRT) together with rational numbers
reconstruction~\cite[\S$5.10$]{GaGe13}.\\

\textit{Implementation of the option ``hybrid'':}
Computing the bounds $(b_t, b_{z_0})$ in step~$\mathbf{1}$
is done using the option~``elimination'' of~Section~\ref{sec:elim}.
More precisely, computing~$b_{z_0}$ is done by specializing~$t$
at a random value~$\theta$ of some random
prime field~$\mathbb{F}_p$ and performing the option ``elimination'' with base field~$\mathbb{F}_p$:
it outputs~$R(\theta, z_0)$ and yields~$b_{z_0}= \deg_{z_0}(R(\theta, z_0))$
(a similar computation is done for computing~$b_t$).
The generation of terms is done using a divide and conquer approach which
computes the first terms of~$F(t, u)\bmod t^{2b_tb_{z_0}+1}$ before
specializing them to~$u=a$.
Finally, the guessing is done by computing Hermite-Padé approximants. For this, we use the function
$\opname{seriestoalgeq}$ of the Maple
package~\href{https://perso.ens-lyon.fr/bruno.salvy/software/the-gfun-package/}{gfun}~\cite{gfun}.

  \section{Getting started}\label{sec:start}
  \subsection{Installation and customization}

  {\em Files provided:
  \href{https://mathexp.eu/notarantonio/}{DDE-Solver}}
  is available in the ``.mla'' format. If one wishes
  to customize the package,
  the Maple scripts ``ddesolver.mpl'' and ``build.mpl''
  can be downloaded on the dedicated
  \href{https://github.com/HNotarantonio/ddesolver}{github webpage}. Used as follows, they allow
  to modify the package by modifying ``ddesolver.mpl'' and to generate
  the corresponding version of~``ddesolver.mla'' by executing~``build.mpl''.\\

  \noindent {\em Customization of \href{https://mathexp.eu/notarantonio/}{DDE-Solver}:}
  Replacing ``PATH/TO'' in the Maple file ``build.mpl'' by the relevant path to the file
  where ``ddesolver.mpl'' and ``build.mpl'' are located allows one, by 
  executing ``build.mpl'', to generate a new version of ``ddesolver.mla''.\\

  \noindent {\em Loading~\href{https://mathexp.eu/notarantonio/}{DDE-Solver}:}
  The Maple variable
  \href{https://fr.maplesoft.com/support/help/maple/view.aspx?path=libname}{libname}
  shall be set so that ``ddesolver.mla'' is located in a visible place.\\

 \code{
  libname :=
  "/home/notarantonio/ddesolver/lib", libname:
}
 \noindent  Once libname has been correctly set up, one executes in Maple\\
  \code{with(ddesolver);} \noindent in order to load and use the package.
  
  \subsection{Using~\ddesolver}

  \paragraph{\textbf{Input/Output syntax.}}

  The arguments of~\ann are $(P, k, \text{var})$,
  where~$P\in\mathbb{Q}[x, z_0, \ldots, z_{k-1}, t, u]$
  is the polynomial in~\eqref{eqn:initpoleqn},~$k$ is the order of the considered DDE and
  var is~$[x, z_0, \ldots, z_{k-1}, t, u]$ (precisely in that order).\\

  The output of~\ann$\!(P, k, \text{var})$
  is some nonzero polynomial~$R\in\Q[t, z_0]$ such
  that~$R(t, F(t, a))=0$.

\begin{tcolorbox}[colback=green!5!white,colframe=green!75!green,]
  \begin{example}\label{ex:3const}
    Consider the DDE~\cite[Eq.~$29$]{BMJ06} of the enumeration of~$3$-constellations:
    \begin{align}\label{eqn:3const}
  \footnotesize    F(t, u) = 1 + tuF(t, u)^3 &+ tu(2F(t, u) + F(t, 1))\frac{F(t, u)-F(t, 1)}{u-1}\\
      \nonumber      &+ tu\frac{F(t, u) - F(t, 1) - (u-1)\partial_uF(t, 1)}{(u-1)^2}.
    \end{align}
    Multiplying~\eqref{eqn:3const} by~$(u-1)^2$ yields
  {\footnotesize  \begin{align*}
     0=(u-1)^2(1-F(t, u) + tuF(t, u)^3) &+ tu(u-1)(2F(t, u) + F(t, 1))(F(t, u)-F(t, 1))\\
      &+  tu(F(t, u) - F(t, 1) - (u-1)\partial_uF(t, 1)).
    \end{align*}
    }
    Thus~$P := (u-1)^2(1-x+tux^3) +tu(u-1)(2x+z_0)(x-z_0)
    +tu(x-z_0-(u-1)z_1)$ and~$k:=2$. We continue the analysis with Maple
    \code{ \raggedright $P := (u-1)^2(1-x+tux^3) +tu(u-1)(2x+z_0)(x-z_0)
    +tu(x-z_0-(u-1)z_1);$\\
{  \raggedright   $with(ddesolver):$}\\
{ \raggedright   $annihilating\_polynomial(P, 2, [x, z_0, z_1, t, u]);$}\vspace{0.1cm}
  { \centering $(16tz_0^2-8tz_0+t-16)(81t^2z_0^3-81t^2z_0^2+27t^2z_0+18tz_0^2-3t^2-66
    tz_0+47t+z_0-1)$}}
    Thus~$R := (16tz_0^2-8tz_0+t-16)(81t^2z_0^3-81t^2z_0^2+27t^2z_0+18tz_0^2-3t^2-66
    tz_0+47t+z_0-1)$ is an annihilating polynomial of~$F(t, a)$.
    A direct analysis on the solutions of~$R(t, z_0)$ that are finite at~$t=0$ shows that
    the second factor of~$R$ is the minimal polynomial of~$F(t,1)$.
    \end{example}
\end{tcolorbox}

  \paragraph{\textbf{Options.}}
  It is possible to benefit in practice
  from two options: the choice of the \emph{algorithm}, and the choice
  of the \emph{variable} (either~$t$ or~$z_0$) on which we perform evaluation--interpolation.\\
  
\noindent \textit{Choice of the algorithm.}
Four algorithms are implemented: ``duplication'', ``elimination'', ``geometry''
(only when~$k=2$) and ``hybrid''.
By default, the algorithm used is ``elimination''.\\

\noindent \textit{Choice of the variable.}
The two choices are~$t$ and~$z_0$. The default choice is~$t$.\\

The choice of the algorithm and of the variable
on which we perform evaluation--interpolation
can be made by
executing~\[\text{\ann}\!(P, k, \text{var}, \text{algorithm}, \text{variable}).\]
Note that these two options must be either not specified at all,
 or specified in the same call.

{\centering\begin{tcolorbox}[colback=green!5!white,colframe=green!75!green,]
      \begin{example} Continuing Example~\ref{ex:3const}
    with~$\opname{algorithm}=$~``geometry'' and~$\opname{variable}=z_0$.
\code{ 
 \raggedright   $annihilating\_polynomial(P, 2, [x, z_0, z_1, t, u], ``geometry'', z_0);$\vspace{0.1cm}
  { \centering $(16tz_0^2-8tz_0+t-16)(81t^2z_0^3-81t^2z_0^2+27t^2z_0+18tz_0^2-3t^2-66
    tz_0+47t+z_0-1)$}}
\end{example}
\end{tcolorbox}
}

We refer to the next section for a discussion on the choice of the parameter~$\operatorname{variable}$.

  \section{Examples}\label{sec:examples}

\subsection{Impact of the option $\opname{variable}$}\label{sec:opt_var}
%
The motivation of this study is the following.
Assume that we perform evaluation--interpolation on~$z_0$. The chosen algorithm
(e.g.~``elimination'' from~Section~\ref{sec:idea_duplication}) will be run
a first time with~$t$ specialized (in order to get the degree in~$z_0$ of the output~$R\in\Q[t, z_0]$),
then it will be run
multiple times at~$z_0$ specialized (in order to evaluate--interpolate~$R$ w.r.t~$z_0$).
Hence if the algorithm is faster at~$t$ specialized than at~$z_0$ specialized (and if the partial
degrees of the output~$R$ are ``close''), then the choice of the evaluation--interpolation variable
has a significant impact on the timings. This is what happens below.\\

{\centering\begin{tcolorbox}[colback=green!5!white,colframe=green!75!green,]
  \begin{example}\label{ex:3-Tamari} Consider the DDE of the enumeration of~$3$-Tamari
    lattices~\cite[Prop.~$8$]{BmFuPr11}:
    \begin{equation}\label{eqn:3tamari}
      \scriptsize   F(t, u) = u +tu F(t, u)\frac{F(t, u)\frac{F(t, u)\frac{F(t, u) - F(t, 1)}{u - 1}
          - F(t, 1)\partial_uF(t, 1)}{u - 1}
    - \frac{F(t, 1)^2\partial_u^2F(t, 1)}{2} - F(t, 1)\partial_uF(t, 1)^2}{u - 1},
    \end{equation}
    Multiplying~\eqref{eqn:3tamari} by~$(u-1)^3$ yields the functional equation
    $P(F(t, u), F(t, 1), \partial_uF(t, 1), t, u)=0,$
    where~$P := -tu^3xz_0^2z_2-2tu^3xz_0z_1^2-2tu^2x^2z_0z_1+2tu^2xz_0^2z_2+4tu^2
xz_0z_1^2+2tux^4-2tux^3z_0+2tux^2z_0z_1-tuxz_0^2z_2-2tuxz_0z_1^2+
2u^4-2u^3x-6u^3+6u^2x+6u^2-6ux-2u+2x\in\Q[x, z_0, z_1, z_2, t, u]$.
We thus have in Maple
\code{ \raggedright $ P := -tu^3xz_0^2z_2-2tu^3xz_0z_1^2-2tu^2x^2z_0z_1+2tu^2xz_0^2z_2+4tu^2
xz_0z_1^2+2tux^4-2tux^3z_0+2tux^2z_0z_1-tuxz_0^2z_2-2tuxz_0z_1^2+
2u^4-2u^3x-6u^3+6u^2x+6u^2-6ux-2u+2x:$\\
{  \raggedright   $with(ddesolver):$}\\
{ \raggedright   $time(annihilating\_polynomial(P, 2, [x, z_0, z_1, z_2, t, u]));$
\\ \hfill \#~$\opname{variable}=t$ (implicitly chosen)}
\\\vspace{0.2cm}
\qquad\qquad\qquad\qquad\qquad\qquad\qquad\;\;\;{\centering $2089.599$}\\
{ \raggedright   $time(annihilating\_polynomial(P, 2, [x, z_0, z_1, z_2, t, u], ``elimination'', z_0));$
                                       \\\hfill   \#~$\opname{variable}=z_0$}\\\vspace{0.2cm}
\qquad\qquad\qquad\qquad\qquad\qquad\qquad\;\;\;{\centering  $141.897$}\\
{\raggedright $R := annihilating\_polynomial(P, 2, [x, z_0, z_1, z_2, t, u],
``elimination'', z_0);$}\\\vspace{0.2cm}
\centering
$R := t^5z_0^{16}+135t^4z_0^{13}+1024t^4z_0^{12}+7290t^3z_0^{10}-1762560t^3z_0^9+393216t^3z_0^8
+196830t^2
z_0^7+111694464t^2z_0^6+580976640t^2z_0^5+t(67108864t+2657205)z_0^4-661978656tz_0^3+4721836032
tz_0^2+(-8371830784t+14348907)z_0+4294967296t-14348907$}
\end{example}
\end{tcolorbox}
}

\subsection{Practical results of the function~\ann}\label{sec:ddesolver_practical_timings}

We provide below a table gathering timings obtained after using the function
\begin{center} \ann\!.\end{center}
Each column is associated with a DDE taken in the literature~\cite{Bernardi08,BMJ06,BmFuPr11,BMC23}.
For each of these DDEs, we precise: its order~$k$, the variable on which we perform
evaluation--interpolation, the algorithm which is used in the call to~\ann\!, and finally the
bi-degree of the output polynomial~$R\in\Q[t, z_0]$.

\begin{landscape}
\small
\vspace*{\fill}
  \begin{center}
  
  \begin{table}
  {\footnotesize
\begin{tabular}{|c|| c c| c c| c c| c c| c c| c c| c c| c c|}\label{table}
  Data
  &  \multicolumn{2}{c|}{$\mathbf{[1]}$} & \multicolumn{2}{c|}{$\mathbf{[2]}$}
  & \multicolumn{2}{c|}{$\mathbf{[3]}$} & \multicolumn{2}{c|}{$\mathbf{[4]}$}
  & \multicolumn{2}{c|}{$\mathbf{[5]}$} & \multicolumn{2}{c|}{$\mathbf{[6]}$}
  & \multicolumn{2}{c|}{$\mathbf{[7]}$} & \multicolumn{2}{c|}{$\mathbf{[8]}$} \\
  \hline
  \hline
  $k$        & \multicolumn{2}{c|}{$2$} & \multicolumn{2}{c|}{$2$} &
  \multicolumn{2}{c|}{$2$} & \multicolumn{2}{c|}{$2$}
  & \multicolumn{2}{c|}{$3$} & \multicolumn{2}{c|}{$3$} & \multicolumn{2}{c|}{$3$}
  & \multicolumn{2}{c|}{$4$} \\
  \hline
  eval--interp variable & $z_0$& $t$ & $z_0$& $t$ &$z_0$& $t$ &$z_0$& $t$ &$z_0$& $t$ &$z_0$& $t$
  &$z_0$& $t$ &$z_0$& $t$\\
  \hline
  ``duplication''       &  $0.5$s    & $0.7$s    & $1.6$s
  & $2.3$s & $41$m$^*$   &  $10$m$^*$   &  $2.7$s   & $1.5$s   &  $2$m$45$s   &  $2$m$50$s & $13$h
  & $27$h  & $\infty$    &  $\infty$   &  $\infty$    &  $\infty$   \\
  
  ``elimination'' &  \color{teal}{$\mathbf{0.2}$\textbf{s}}    &  $0.3$s
  &  \color{teal}{$\mathbf{0.4}$\textbf{s}}    &  \color{teal}{$\mathbf{0.4}$\textbf{s}}
  & $1$h$20$m$^*$   &   $4$m$^{*}$  &   $1.7$s  &
  $1.2$s
  &   $1$m$10$s
  &  \color{teal}{$\mathbf{45}$\textbf{s}}  &
  \color{teal}{$\mathbf{2}$\textbf{m}$\mathbf{20}$\textbf{s}}
  &   $35$m      &    $1$m       &  $7$m$30$s
  &     $2$d$19$h     &     $2$d     \\
  
  ``geometry''     &  $2$s   &  $2.2$s    & $1.6$s    &   $0.9$s    &  $54$m$^*$   &
  \color{teal}{$\mathbf{2}$\textbf{m}}$^*$  & $0.8$s &
  \color{teal}{$\mathbf{0.7}$\textbf{s}}
  & $\times$  & $\times$    & $\times$ &$\times$ &$\times$ &$\times$ &
  $\times$ &$\times$\\
  ``hybrid'' &   \multicolumn{2}{c|}{$50$s}  &  \multicolumn{2}{c|}{$2$m} & \multicolumn{2}{c|}{$\infty$}
  &\multicolumn{2}{c|}{$1.4$s} & 
  \multicolumn{2}{c|}{$18$s} & \multicolumn{2}{c|}{$1$h$42$m} &
  \multicolumn{2}{c|}{\color{teal}{$\mathbf{34}$\textbf{s}}} &
  \multicolumn{2}{c|}{\textcolor{teal}{$\mathbf{2}$\textbf{h}$\mathbf{41}$\textbf{m}}}\\
  \hline
  ($\deg_{t}(R)$, $\deg_{z_0}(R)$) &  \multicolumn{2}{c|}{($4,9$)}  &  \multicolumn{2}{c|}{($2,3$)}
  & \multicolumn{2}{c|}{($132, 6$)}
  &\multicolumn{2}{c|}{($3, 5$)} & 
  \multicolumn{2}{c|}{($3, 7$)} & \multicolumn{2}{c|}{($5, 16$)} &
  \multicolumn{2}{c|}{($2, 4$)} & \multicolumn{2}{c|}{($3, 9$)}\\
  \hline\end{tabular}
}    \end{table}\vspace{-1cm}

 {\centering \textbf{Table~$1$:} Practical results\footnote{
 All computations are conducted using Maple on a computer equipped
  with Intel® Xeon® Gold CPU 6246R v4 @ 3.40GHz and 1.5TB~of~RAM~with~$1$~thread.}
    \footnote{Gr\"obner bases computations are performed
using the C library~\href{https://msolve.lip6.fr/}{msolve}~\cite{msolve}, implemented 
by Berthomieu, Eder and Safey El Din.}
    of\;\ann on DDEs from the
    literature.}

  \begin{itemize}
\item $\mathbf{[1]}$: Enumeration of $2$-Tamari lattices, \cite[Prop.~$8$]{BmFuPr11} for~$m=3$,
\item $\mathbf{[2]}$: Enumeration of $2$-greedy Tamari intervals,
  \cite[Prop.~$3.1$]{BMC23} for~$m=2$,
\item $\mathbf{[3]}$: Enumeration of non-separable near-triangulations in which all intern
  vertices have degree at least~$5$, \cite[Prop. $4.3$, Eq. $22$]{Bernardi08},
\item $\mathbf{[4]}$: Enumeration of~$3$-constellations, \cite[Equation~$29$]{BMJ06},
\item $\mathbf{[5]}$: Enumeration of~$4$-constellations, \cite[Prop.~$12$]{BMJ06} for~$m=4$,
\item $\mathbf{[6]}$: Enumeration of~$3$-Tamari lattices, \cite[Prop.~$8$]{BmFuPr11} for~$m=4$,
\item $\mathbf{[7]}$: Enumeration of~$3$-greedy Tamari intervals,
    \cite[Prop.~$3.1$]{BMC23} for~$m=3$,
  \item $\mathbf{[8]}$: Enumeration of~$5$-constellations.
    \cite[Prop.~$12$]{BMJ06} for~$m=5$.
    \end{itemize}
\end{center}
\vspace{0.4cm}
\begin{itemize}
\item $\infty$: The computations did not finish within~$5$ days,
\item $\times$: The algorithm is not implemented for~$k>2$,
\item $\cdot^*$: We added the condition $m'\cdot(1+t^3)(1-t^3)t-1=0$ to the input polynomial constraints.
\end{itemize}
\vspace{0.4cm}

\end{landscape}

From Table~$1$, we strengthen the message of~Section~\ref{sec:opt_var} and 
draw some conclusions:

\begin{itemize}
\item {\em Choice of the algorithm}
  The best algorithm to choose depends in an important way on the studied DDE. Even
among DDEs of the same family (e.g. DDEs~$[\mathbf{4}], [\mathbf{5}]$ and~$[\mathbf{8}]$)
the answer is not clear: the first two
are solved faster with the option $\opname{algorithm}=$``elimination'' while the last one is solved
faster with the option~$\opname{algorithm}$= ''hybrid''.
The efficiency of this last option~$\opname{algorithm}=$''hybrid'' comes from the fact that
the output polynomial has small bidegree $(3, 9)$, thus allowing a fast computation of the first terms
of~$F(t, 1)$; and at the same time, that the computations are hard with the other algorithms.

This impossibility to predict which algorithm is more efficient for an input DDE
justifies that we let the user choose the algorithm that shall be applied.\vspace{0.2cm}
\item {\em Choice of the variable} 
  As in~Section~\ref{sec:opt_var}, the table shows
  that the choice of the variable ($z_0$ or~$t$)
  on which we perform evaluation--interpolation is important.
  Let us pick one example to illustrate it. For the DDE~$[\mathbf{3}]$,
  performing evaluation--interpolation over~$t$ with~$\opname{algorithm}=$''geometry''
  takes~$2$ minutes while it takes~$54$ minutes when performing evaluation--interpolation over~$z_0$.
  With the notations of~Section~\ref{sec:geometry}, the reason is that
  at~$z_0$ specialized to some random~$\theta\in\mathbb{F}_p$,
  the main computational cost (in addition to the Gr\"obner basis computation)
  is to compute the discriminant w.r.t.~$T$ of~$(\chi_{z_1})|_{z_0=\theta}\in\Q[t, T]$,
  whose partial degrees are~$6$ in~$T$ and~$423$ in~$t$.
  In the same time at~$t$ specialized to some random~$\nu\in\mathbb{F}_p$,
  the polynomial~$(\chi_{z_1})|_{t=\nu}\in\Q[z_0, T]$
  has degree~$15$ in~$z_0$ and~$6$ in~$T$: it yields a faster computation of its
  discriminant w.r.t. $T$.

  As all these intermediate computational data are hard to predict
  from the input DDE, we let the user choose on which
  variable we shall perform evaluation--interpolation.\vspace{0.2cm}
\item
{\em Solving DDEs previously out of reach!} The DDEs~$[\mathbf{7}]$,~$[\mathbf{8}]$
are solved via theoretical arguments in (resp.)~\cite{BMC23,BMS00}. Here, it is the first time
that they are solved in an automatized~way!
\end{itemize}

\section{Appendix}

\subsection{Solving~$\mathbf{3}$-constellations using~Section~\ref{sec:idea_duplication}}\label{annex:dup}

  The DDE~\cite[Eq.~$29$]{BMJ06} associated with the enumeration of~$3$-constellations
  is given by
  \begin{align}\label{eqn:3constbis}
    F(t, u) = 1 + tuF(t, u)^3 +& tu(2F(t, u)+F(t, 1))\frac{F(t, u)-F(t, 1)}{u-1}\\
  \nonumber  +&tu\frac{F(t, u)-F(t, 1)-(u-1)\partial_uF(t, 1)}{(u-1)^2}
    \end{align}
  Multiplying~\eqref{eqn:3constbis} by~$(u-1)^2$ yields
$P(F(t, u), F(t, 1), \partial_uF(t, 1), t, u)=0$,
  where~
  {\footnotesize\[P := (u-1)^2(1-x+tux^3) +tu(u-1)(2x+z_0)(x-z_0)
    +tu(x-z_0-(u-1)z_1)\in\Q[x, z_0, z_1, t, u].\]
    }

   \noindent Define the polynomial system
    \begin{align*}
\mathcal{S}_{\textsf{dup}}:=  (&P(x_1, z_0, z_1, t, u_1)=0, \partial_xP(x_1, z_0, z_1, t, u_1)=0,
  \partial_uP(x_1, z_0, z_1, t, u_1)=0, \\
  &P(x_2, z_0, z_1, t, u_2)=0,
\partial_xP(x_2, z_0, z_1, t, u_2)=0, \partial_uP(x_2, z_0, z_1, t, u_2)=0,\\
&m\cdot (u_1-u_2)\cdot(u_1-1)\cdot(u_2-1)\cdot u_1\cdot u_2-1=0),
\end{align*}

 \noindent    and set~$\mathcal{I}_\textsf{dup}\subset\Q(t)[m, x_1, x_2, z_0, z_1, u_1, u_2]$
    to be the ideal generated by the polynomials~in~$\mathcal{S}_{\textsf{dup}}$.\\
     
    The algorithm described in~Section~\ref{sec:idea_duplication}:
    \begin{enumerate}
\item Computes a generator~$R\in \mathcal{I}_\textsf{dup}
  \cap\Q(t)[z_0]$ by applying~Fact~\ref{fact:fglm} and finds
  \[R = (81t^2z_0^3-81t^2z_0^2+27t^2z_0+18tz_0^2-3t^2-66tz_0+47t+z_0-1)
  (16tz_0^2-8tz_0+t-16).\]
\end{enumerate}

\subsection{Solving~$\mathbf{3}$-constellations using~Section~\ref{sec:elim}}\label{annex:elim}

Define~$\mathcal{S} := (P, \partial_xP, \partial_uP, m\cdot u(u-1)-1)$, for~$P$
as in~\ref{annex:dup}. Denote by~$\mathcal{X}\subset\overline{\Q(t)}^5$ the solution set
of the conditions~$P=0 \wedge \partial_xP =0 \wedge \partial_u P=0 \wedge u(u-1)\neq0$
in~$\overline{\Q(t)}^4$.\\

\textbf{Goal:} Characterize~$\mathcal{F}_k(u, \mathcal{X}) := \{ \bm{\alpha}\in\overline{\Q(t)}^2 \;|\;
\#_u(\mathcal{X}, \bm{\alpha})\geq 2\}$ as the solution set of some polynomial constraints and
deduce a nonzero~$R\in\Q[t, z_0]$ such that
\[R(t, F(t, 1))=0.\]

The algorithm described in~Section~\ref{sec:elim}:

\begin{enumerate}
\item Characterizes the projection of~$\mathcal{X}$
  onto the~$(u, z_0, z_1)$-coordinate space:
  \code{
    \raggedright $G$ := Groebner[Basis]($\mathcal{S}$, lexdeg($[m, x], [u, z_0, z_1, t]$)):\\
    \raggedright $E$ := Groebner[Basis](remove(has, $G$, $\{m, x\}$), lexdeg($[u], [z_0, z_1, t]$)):\\
    \raggedright $Leading$ := select(has, $G$, $\{m, x\}$):}
  The first element of~$Leading$ has
  the form
  { 
  \[1534(u-1)\cdot x + (\text{polynomial in }u, z_0, z_1, t).\]
  }
  As
    $\mathcal{X}\cap \{u=1\}=\emptyset$, it results from the extension theorem that
    the zero set of~$E$ in~$\overline{\Q(t)}^3$ is precisely the projection of~$\mathcal{X}$
    onto the~$(u, z_0, z_1)$-coordinate space.
  \item Identifies algebraic conditions in~$E$ describing the cardinality of the fiber
    in~$\mathcal{F}_k(u, \mathcal{X})$:
    \code{\raggedright\textbf{\# Dealing with the polynomials of degree~$\mathbf{< 2}$
        in~$\mathbf{u}$ in~$E$:}\\
      \raggedright $NewConditions$ := [];\\
      \raggedright for $i$ to nops($E$) do \\
      \raggedright \qquad if degree($E$[$i$], $u$) $< 2$ then \\
      \raggedright \qquad \qquad $NewConditions$ :=
                   [op($NewConditions$), coeffs($E$[$i$], $u$)]:\\
      \raggedright \qquad fi:\\
      \raggedright od:\\
      \raggedright $E$ := Groebner[Basis]([op($E$), op($NewConditions$)],
      lexdeg([$u$], [$z_0, z_1, t$])):
      \\\vspace{0.5cm}
      \raggedright\textbf{\# Dealing with the polynomials of degree~$\mathbf{\geq 2}$ in~$\mathbf{u}$
        in~$E$:}\\
      \raggedright map($p\rightarrow$degree($p, u$), $E$); \# prints the degrees in~$u$ of
      the polynomials in~$E$\\\vspace{0.2cm}
      \centering   $[0, 0, 0, 0, 0, 0, 0, 0, 0, 2, 2, 2, 2, 2, 2, 2, 2, 2, 2, 3, 3, 3]$\\
      \raggedright \#
      In order to have $\#_u(\mathcal{X}, \bm{\alpha})\geq 2$
      for a given~$\bm{\alpha}\in\overline{\Q(t)}^2$,\\
      \raggedright \# we shall exclude the solution set
      of~$\opname{disc}_u(E[10])$:\\
      \raggedright  $sat$ := discrim($E[10]$, $u$):\\
    }
    Denote by~$c_{10}(z_0, z_1)\in\Q(t)[z_0, z_1]$ the coefficients of~$u^2$ in~$E[10]$.
    As~$V(E)\cap \{c_{10} = 0\} = \emptyset$, the algorithm deduces from
    Theorem~\ref{thm:extension} that
    $\mathcal{F}_k(u, \mathcal{X})$ is the solution set of the polynomial equations given by~$E$ and
    of the inequation~$sat\neq 0$.
  \item Turns the inequation~$sat\neq0$
    into the equation~$m\cdot sat-1=0$ by introducing an extra variable~$m$.
    Finally, it deduces a nonzero polynomial~$R\in\Q[t, z_0]$ annihilating~$F(t, 1)$.
    \code{
         \raggedright $H$ := Groebner[Basis]([op($E$), $m\cdot sat-1$],
      lexdeg([$m, u, z_1$], [$z_0, t$])):\\
      \raggedright op(remove(has, $H$, $\{m, u, z_1\}$));\\\vspace{0.2cm}
      \centering
            $ (16tz_0^2-8tz_0+t-16)(81t^2z_0^3-81t^2z_0^2+27t^2z_0+18tz_0^2-3t^2-66tz_0+47t+z_0-1)$
    }
    Thus~\[R = (16tz_0^2-8tz_0+t-16)(81t^2z_0^3-81t^2z_0^2+27t^2z_0+18tz_0^2-3t^2-66tz_0+47t+z_0-1).\]
  \end{enumerate}

\subsection{Solving~$\mathbf{3}$-constellations using~Section~\ref{sec:geometry}}\label{annex:geom}

Denote by~$P\in\Q[x, z_0, z_1, t, u]$ the polynomial in~\ref{annex:dup}
and by~$\mathcal{S}$ the set of polynomials~$(P, \partial_x P, \partial_u P, m\cdot tu(u-1)-1)
\subset\Q(t)[m, x, z_0, z_1, u]$. Moreover, recall that~$\mathcal{X}\subset\overline{\Q(t)}^4$ is the
solution set of the constraints~$P=0\wedge \partial_xP=0\wedge \partial_uP=0\wedge u(u-1)\neq0$ and
 that

{\footnotesize\[\mathcal{S}_k(\mathcal{X}) := \{\bm{\alpha} = (\alpha_0, \ldots, \alpha_{k-1})
\in\algebraicclosure^k \;\,|\,\;
\bm{\alpha}\in\pi(\mathcal{X}) \;\,\wedge\;\, \mathcal{X} \cap\;
\pi_{\check{z_1}}^{-1}((\alpha_0, \alpha_2, \ldots, \alpha_{k-1}))\geq
k\}.\]
}
%

\textbf{Goal:} Characterize~$\mathcal{S}_k(\mathcal{X})$ as the solution set of some
polynomial constraints and deduce a nonzero~$R\in\Q[t, z_0]$ such that~$R(t, F(t, 1))=0$.\\

The algorithm described in~Section~\ref{sec:geometry}:
\begin{enumerate}
\item Computes a Gr\"obner basis of the ideal generated by~$\mathcal{S}$ for the
  order~$\succ_{\text{bgrevlex}}$:
\code{
 \raggedright   $G$ := Groebner[Basis]($\mathcal{S}$, lexdeg([$m, x, u, z_1$], [$t, z_0$])):
}
\item Computes the matrix~$M_{z_1}$ of the multiplication map~$m_{z_1}:f\mapsto z_1\cdot f$ in the quotient
  ring $\Q(t, z_0)[m, x, u, z_1]/\langle j(\mathcal{S})\rangle$, where~$j$ is the usual inclusion map
  of~$\Q(t)[m, x, z_0, z_1]$ in~$\Q(t, z_0)[m, x, u, z_1]$:
  \code{
{\raggedright {\footnotesize $L_1, L_2$ := Groebner[NormalSet](subs($t = rand()$, $z_0 = rand()$,
  $G$), tdeg($m, x, u, z_1$)):}}\\
\raggedright  $M$ := Groebner[MultiplicationMatrix]($z_1$, $L_1$, $L_2$,
  $G$, tdeg($m, x, u, z_1$)):
  }
\item Computes the characteristic polynomial of~$M_{z_1}$ and defines~$\chi_{z_1}$ to be its numerator:
  \code{
{ \raggedright $\xi_{z_1}$ := \text{LinearAlgebra[CharacteristicPolynomial]}(M, T):}\\
\raggedright  $\chi_{z_1}$ := \text{factor}(\text{numer}($\xi$));\\\vspace{0.2cm}
%
{\centering$(\mb{T}t+tz_0^2-1)\cdot (729\mb{T}^4t^4+108t^3(27tz_0^2-9tz_0+2t-12)\mb{T}^3+2t(2187t^3z_0^4-1350t^3z_0^3+459t^3z_0^2-72t^3z_0-2376t^2z_0^2+8t^3+1728t^2z_0-492t^2-36tz_0+276t-8)\mb{T}^2+4t(729t^3z_0^6-621t^3z_0^5+261t^3z_0^4-55t^3z_0^3-1404t^2z_0^4+6t^3z_0^2+1536t^2z_0^3-560t^2z_0^2-36tz_0^3+68t^2z_0+640tz_0^2-8t^2-916tz_0-8z_0^2+328t+8z_0)\mb{T}+ 16+358t^4z_0^6-84t^4z_0^5+9t^4z_0^4-1160t^3z_0^4-3224t^2z_0^3+1320t^2z_0^2-128t^2z_0-720tz_0^2+729t^4z_0^8-756t^4z_0^7-2160t^3z_0^6+2752t^3z_0^5-72t^2z_0^5+2104t^2z_0^4-16tz_0^4+16t^2-32z_0-544t+96tz_0^3-24t^3z_0^2+208t^3z_0^3+1184tz_0+16z_0^2)$
  } }
  \item
    Thus~$\mathcal{S}_k(\mathcal{X})$ is the solution set of the constraints
    $P=0\wedge \partial_xP=0 \wedge \partial_u P=0 \wedge u(u-1)\neq0\wedge
    \opname{disc}_T(\chi_{z_1})=0$.
    An annihilating polynomial of~$F(t, 1)$ is given by
    the discriminant~$R\in\Q[t, z_0]$ of~$\chi_{z_1}$ with respect to~$T$, that is
\begin{align*}
\tiny  R =&\; 12230590464t^{14}\cdot
  (16tz_0^2-8tz_0+t-16)^2\cdot(1+(3z_0-1)^2t^2+(-6z_0-14)t)^3\\
&\cdot (tz_0+1)^{12}\cdot(81t^2z_0^3-9t(9t-2)z_0^2+(27t^2-66t+1)z_0-3t^2+47t-1)^2.
  \end{align*}

\end{enumerate}

\subsection{Solving~$\mathbf{3}$-constellations using~Section~\ref{sec:hybrid}}\label{annex:hybrid}

Consider~\eqref{eqn:3constbis} and the polynomial~$P\in\Q[x, z_0, z_1, t, u]$ of~\ref{annex:dup}.
Define the set of polynomials $\mathcal{S} := (P, \partial_xP, \partial_uP, m\cdot u(u-1)-1)\in
\Q(t)[m, x, z_0, z_1, u]$.\\

\textbf{Goal:} Compute a nonzero polynomial~$R\in\Q[t, z_0]$ such that~$R(t, F(t, 1))=0$.\\

The algorithm described in~Section~\ref{sec:hybrid}:
    
    \begin{enumerate}
    \item
    Applies twice the strategy described in~\ref{annex:elim} to~$P$. It draws
    at random a prime number~$p$ and some~$\theta\in\mathbb{F}_p$, say~$p = 12301$ and~$\theta = 1328$,
    and applies~\ref{annex:elim} with~$\mathcal{S}$ replaced by~$(\mathcal{S})|_{t = \theta}$
    and~$\Q$ replaced by~$\mathbb{F}_p$. It outputs
    \[
    7957(z_0+11829)(z_0^3+4863z_0^2+8711z_0+3012)(z_0+6622)\]
    \noindent The degree being~$5$, it sets~$b_{z_0} = 5$.
    Similarly, it obtains~$b_t = 3$.
      
      \item  Computes the~$32$ first terms of~$F(t, 1)$:
       \[F_1 = 1 + t + 6t^2 + 54t^3 + 594t^4 + \cdots + 913075994651156584840651326232625946t^{31}
       \bmod t^{32}\]
            
    \item Uses the function~$\opname{seriestoalgeq}$ of the Maple package
      \href{http://perso.ens-lyon.fr/bruno.salvy/software/the-gfun-package/}{gfun}
      \cite{gfun}:
      \code{
        \raggedright{\scriptsize$\opname{collect}(\opname{subs}(T(t) = z_0,
        \opname{gfun}[\opname{seriestoalgeq}](\opname{series}(F_1, t,
        \opname{degree}(F_1, t)),        T(t)))[1], z_0, \opname{factor});$}\\\vspace{0.2cm}
        \centering  $81t^2z_0^3-9t(9t-2)z_0^2+(27t^2-66t+1)z_0-3t^2+47t-1$
      }
    \noindent  The above polynomial, denoted~$M$,
    is a candidate for annihilating~$F(t, 1)$.
    
      \item Proves the guessed polynomial with the below lines
      \code{\raggedright series(subs($z_0 = F_1, M$), $t$, $31$);
        \# Compute $M(t, F(t, 1) \bmod t^{32})$\\\vspace{0.2cm}
        \centering  $O(t^{31})$}
      Thus~$R = 81t^2z_0^3-9t(9t-2)z_0^2+(27t^2-66t+1)z_0-3t^2+47t-1$.\\
\end{enumerate}

\paragraph{\textbf{Acknowledgment.}} 
The author warmly thanks Alin Bostan and Mohab Safey El Din for various
important suggestions that yielded the current version of this paper. The author was
supported by the ANR-19-CE40-0018 project De Rerum Natura, the French–Austrian
ANR-22-CE91-0007 FWF I6130-N project EAGLES and the ANR-21-CE48-0007 project IsOMa. 
Also, the author acknowledge support of the Institut Henri Poincaré (UAR 839 CNRS-
Sorbonne Université), and LabEx CARMIN (ANR-10-LABX-59-01).
The author thanks Mireille Bousquet-Mélou for her suggestions.\!

\end{document}